\input amstex

\documentstyle{amsppt}
%\langno0   %%% 0 - English, 1 - Ukrainian, 2 - Russian, 3 - Francais

\pageno=99
%{\sixrm ""Š 515.12}
%\vs
\NoBlackBoxes

\def\vs{\vskip3pt}
\def\qed{\hfill$\square$}

\def\U{{\Cal U}}
\def\Qi{Q^\infty}

\def\C{{\Cal C}}
\def\conv{\operatorname{conv}}
\def\IN{{\Bbb N}}
\def\IR{{\Bbb R}}
\def\vsk{\vskip10pt}
\def\bs{\setminus}
\def\Ri{\IR^\infty}
\def\spa{\operatorname{span}}

\def\dlim{\varinjlim}
\def\Cfd{\C^\infty _{fd}}

\def\emptyset{\varnothing}

\topmatter
\title
On linear topological  spaces (linearly) homeomorphic to $\IR^\infty$
\endtitle
\rightheadtext{On linear spaces homeomorphic to $\IR^\infty$}
\author
Taras Banakh
\endauthor
\abstract\nofrills{T.\ Banakh. {\it 
On linear topological spaces (linearly) homeomorphic to
$\IR^\infty$}, Matem. Studii. 9:1 (1998), 99--101.
\vskip3pt

We prove that every infinite-dimensional (locally convex) linear topological 
space that can be expressed as a direct limit of finite-dimensional
metrizable compacta is (linearly) homeomorphic to the space
$\IR^\infty=\dlim \IR^n$. 
}\endabstract
\endtopmatter

Given an increasing sequence of topological spaces
$$
X_1\subset X_2\subset\dots \subset X_n\subset\cdots
$$
we define the direct limit topology on $X=\bigcup_{n\in\IN}X_n$ letting
$U\subset X$ to be open if and only if $U\cap X_n$ is open for every
$n\in\IN$.

By $\IR^\infty$ we denote the direct limit of the sequence 
$$
\IR^1\subset \IR^2\subset \IR^3\subset\cdots,
$$
where the embeddings $\IR^n\subset\IR^{n+1}$ are defined by 
$(x_1,\dots,x_n)\mapsto (x_1,\dots,x_n,0)$.

The space $\IR^\infty$ considered with the natural linear operations is a
locally convex linear topological space. Denote by $\Cfd$  the
 class of spaces that are direct limits of
finite-dimensional metrizable compacta. 

In this note we prove the following

\proclaim {Theorem} Any infinite-dimensional (locally convex) linear
topological space $X\in\Cfd$ is (linearly) homeomorphic to $\Ri$.
\endproclaim 

\remark{Remark} A 
metric counterpart of the space $\Ri$ is the linear subspace
$l_2^f=\{(t_i)_{i=1}^\infty \in l_2\mid t_i=0$ for all but finitely many
$i\}$ in the separable Hilbert space $l_2$. 
A part of Theorem (that dealing with homeomorphisms) has its metric
analog:
every infinite-dimensional linear
metrizable space that can be represented as a countable union of
finite-dimensional compacta is homeomorphic to the pre-Hilbert space
$l_2^f$ \cite{CDM}.
In the meantime, the other part of 
Theorem is specific for direct limit
topologies and admits no generalization onto metric locally convex spaces: 
the~linear $\spa([0,1])$ of a~linearly independent arc $[0,1]$ in~$l_2$
is a~locally convex linear metric space (even pre-Hilbert space) which is
a~countable union of finite-dimensional compacta but is not isomorphic to
the~space~$l_2^f$.
\endremark

Theorem results from the~following a~little bit more general results.

\proclaim{Proposition 1}Any convex set $X\in\Cal C_{fd}^\infty$ in 
a~locally convex linear topological space~$L$ is affinely homeomorphic to a~convex
set in~$\IR^\infty$. 
Moreover, if $X=L$ then $X$ is linearly homeomorphic to 
$\IR^\infty$.
\endproclaim

\demo{Proof} If $X$ is finite-dimensional then the statement is trivial,
so we assume that $X$ is infinite-dimensional. Write $X=\dlim X_n$, where
each $X_n$ is a finite-dimensional compactum. Without loss of generality,
$0\in X_1\subset X_2\subset\dots \subset X$.

\proclaim{Claim} For every compactum $C\subset X$ the convex hull
$\conv(C)\subset X$ is finite-dimensional.
\endproclaim

\demo{Proof} Assume on the contrary that $\conv(C)$ is
infinite-dimensional. Then for every $n\in\IN$ there exists
$x_n\in\frac1n\conv(C)\bs X_n$ (recall that $X_n$'s are
finite-dimensional). It follows from the definition of the direct limit
topology on $X$ that the set $A=\{x_n\mid n\in\IN\}$ is closed in~$X$.
Moreover, $0\notin A$. Since the space $L$ is locally convex, there is a
convex neighborhood $U\subset L$ of the origin such that $U\cap
A=\emptyset$. Since the set $C$ is compact, $\frac1nC\subset U$ for some
$n\in\IN$. By convexity of $U$, $x_n\in\conv(\frac1nC)\subset U$. This
contradicts to $x_n\in A$ and $U\cap A=\emptyset$.\hfill$\square$
\enddemo

Applying Claim to the compacta $X_n$'s, one can construct a sequence
$\{x_n\}_{n\in\IN}\subset X$ of linearly independent vectors and a number
sequence $\{m(n)\}_{n\in\IN}$ such that $X_n\subset
L_{m(n)}=\spa\{x_1,\dots, x_{m(n)}\}$ for every $n\in\IN$. Then the
identity map $i\colon X\to \dlim L_{m(n)}$ is continuous. Now remark that for
every $n\in\IN$ the intersection $X\cap L_{m(n)}$ is a metrizable direct
limit of compacta, and hence, $X\cap L_{m(n)}$ is locally compact. Using
this fact, show that $i$ is a topological embedding, and notice that
$\dlim L_{m(n)}$ is affinely homeomorphic to $\IR^\infty$.

If $X=L$ then $X=\bigcup_{n\in\IN}L_{m(n)}$, and consequently, 
$i(X)=\dlim L_{m(n)}$.
\qed
\enddemo

\proclaim{Proposition 2} Every infinite-dimensional convex set $X\in\Cfd$
in a linear topological space $L\in\Cfd$ is homeomorphic to $\IR^\infty$.
\endproclaim

To prove this proposition we need

\proclaim {Lemma} Let $X$ be a convex infinite-dimensional subset of a
linear topological space $L\in\Cfd$. For every compact set $C\subset X$
there is an embedding $e\colon C\times [0,1]\to X$ such that $e(c,0)=c$ for
every $c\in C$.
\endproclaim

\demo {Proof} We will search the embedding $e$ among the maps of the type
$e(c,t)=(1-\frac t2)c+\frac t2x_0$, where $x_0\in X$.

So, let $C\subset X$ be a compactum. Then the set
$D=[0,1]\cdot([1/2,1]C-[1/2,1]C)=\{t(\tau\, c-\tau'c')\mid t\in
[0,1],\;\tau,\tau'\in[1/2,1],\;c,c'\in C\}\subset L$ is also compact.
Moreover, since $L\in\Cfd$, the compactum $D$ is finite-dimensional. Let
$n=\dim(D)$. Since the convex set $X$ is infinite-dimensional it contains
$n+2$ linearly independent vectors $x_1,\dots,x_{n+2}$. Let
$\Delta=\conv\{x_1,\dots,x_{n+2}\}\subset X$ be their convex hull. We claim
that there is $x_0\in\Delta$ such that $x_0\cdot(0,1]\cap D=\emptyset$. 
Assuming the converse, we obtain that
 $\Delta=\bigcup_{k=1}^\infty \Delta_k$, where
$\Delta_k=\{x\in\Delta\mid x\cdot[0,1/k]\subset D\}$
 (remark that
the set $D$ together with a point $d\in D$ contains the interval
$[0,1]d$~).
 Using the
compactness of $D$ prove that each set~$\Delta _k$ is closed in $\Delta$.
Consequently, by the Baire Category Theorem, one of the sets
$\Delta_k$'s (to say $\Delta_m$) is somewhere dense in $\Delta$. Then $\dim
(\Delta_m)=\dim(\Delta)=n+1$ and the map $f(x)=\frac 1m x$ for
$x\in\Delta_m$ determines the embedding $f\colon \Delta_m\to D$ of the
$(n+1)$-dimensional compactum $\Delta_m$ into the $n$-dimensional space
$D$. The obtained contradiction shows that there is $x_0\in \Delta\subset
X$ such that $x_0\cdot (0,1]\cap D=\emptyset$.

Let us show that the map $e\colon C\times [0,1]\to X$ defined by
$e(x,t)=(1-\frac t2)c+\frac t2x_0$ for $(c,t)\in C\times [0,1]$ is an
embedding.
 Since the space $C$ is compact, it suffices to prove that the map
$e$ is injective. Let $(c,t)$ and $(c',t')$ be distinct points in $C\times
[0,1]$. If $t=t'$ then $c\not= c'$ and obviously $e(c,t)\not=e(c',t')$. If
$t\not=t'$ then $\frac {t-t'}2x_0\notin D=[0,1]([1/2,1]C-[1/2,1]C)\ni
(1-\frac {t'}2)c'-(1-\frac t2)c$. Consequently, $\frac {t-t'}2x_0\ne
(1-\frac {t'}2)c'-(1-\frac t2)c$. This yields 
$(1-\frac t2)c+\frac t2x_0\ne(1-\frac {t'}2)c+\frac {t'}2x_0$, i.e., the
map $e$ is injective. \qed
\enddemo 

{\noindent\it Proof of Proposition 2.} Let $X\in\Cfd$ be a convex
infinite-dimensional subset in a linear topological space $L\in\Cfd$. 
To prove that  the space $X$ is homeomorphic to $\Ri$ we will apply 
Sakai's Characterizing Theorem \cite{Sa}. 
According to this theorem, a space $X\in \Cfd$ is homeomorphic to $\IR^\infty$ 
if and only if for every finite-dimensional compact pair $B\subset A$ any
embedding $i\colon B\to X$ extends to an embedding
$\bar i\colon A\to X$.

Fix a finite-dimensional compactum $A$, a closed subset $B$, and an
embedding $i\colon B\to X$.
At first we shall construct a map $\tilde i\colon A\to X$ extending the
embedding $i$. For this fix any metric $d$ on the compactum $A$ and
consider the cover $\U=\{O_d(a,d(a,B)/3)\mid a\in A\bs  B\}$
of the space
$A\bs B$. For every $U=O(a,d(a,B)/3)\in\U$ fix a point $b_a\in B$ with
$d(a,b_a)=d(a,B)$. Let $\{\lambda _U\colon A\bs B\to [0,1]\}_{U\in\U}$ be a
partition of unity of order $\le \dim(A)+1$ inscribed into the cover $\U$.
It can be easily verified that the map $\tilde i\colon A\to X$ defined by
$$
\tilde i(a)=\cases i(a), & a\in B\\
\sum_{U\in \U}\lambda _U(a)i(b_a),& a\in A\bs B\endcases
$$
is the required extension of the embedding $i\colon B\to X$.

Let $j\colon A\bs B\to [0,1]^m$ be an embedding of the quotient space $A/B$ such
that $j(\{B\})=\bold 0=(0,\dots,0)$. It follows from Lemma that there is
an  embedding $e\colon \tilde i(A)\times [0,1]^m\to X$ such that $e(x,0)=x$ for
$x\in\bar i(A)$. Denote by $\pi\colon A\to A/B$ the quotient map. It is easily
seen that the map $\bar i=e(\tilde i,j\circ \pi)\colon A\to X$ defined by $\bar
i(a)=e(\bar i(a),j\circ\pi(a))$, $a\in A$, is the required embedding
extending the embedding $i\colon B\to X$. By \cite{Sa},
 the space $X$ is homeomorphic to $\Ri$.\qed
\vsk

\centerline{REFERENCES}
\vs

\vsk

Department of Mathematics, Lviv University,  

Universytetska 1, Lviv,
290602, Ukraine
\vs

\rightline{Received 15.01.1996}

\rightline{Revised 15.01.1997}
\bye